\documentclass[]{amsart}
\usepackage{amsmath,amsfonts,amssymb,amsthm}
\usepackage{mystyle}  

\usepackage{dsfont}
\def\k{{\mathds k}}
\def\1{{\mathds 1}}
\def\Z{{\mathds Z}}
\def\Q{{\mathds Q}}
\def\C{{\mathds C}}

\renewcommand\^[1]{\widehat{#1}}

\def\form#1#2{\langle#1\,|\,#2\rangle}
\def\bform#1#2{\big\langle#1\,\big|\,#2\big\rangle}
\def\ord{\op{ord}}

\def\mylabel#1{\label {#1}} 

\begin{document}

\bibliographystyle{plain}

\author{Michael Roitman}

\address{University of Michigan\\
Department of Mathematics\\
Ann Arbor, MI 48109-1109}
\email{roitman@umich.edu}

\title{Invariant bilinear forms on a vertex algebra and the radical}

\date\today




\maketitle
\section*{Introduction}
Invariant bilinear forms on vertex algebras have been around for quite
some time now. They were mentioned by Borcherds in \cite{bor} and were
used in many early works on vertex algebras, especially in relation
with the vertex algebras associated with lattices
\cite{dong,flm,kac2}. The first systematic study of invariant forms on
vertex algebras is due to Frenkel, Huang and Lepowsky \cite{fhl}. This
theory was developed further by Li \cite{liform}.

However, these authors imposed certain assumptions on their vertex
algebras which are too restrictive for the applications we have in
mind. Specifically, this paper is motivated by the study of vertex
algebras of OZ type generated by their Griess subalgebras
\cite{gnavoa,V2}. We show that all the results of Li \cite{liform}
hold in the greatest possible generality, basically, as long as the
definitions make sense. We construct a linear space that
parameterizes all bilinear forms on a given vertex algebra with values
in a arbitrary linear space, and also
prove that every invariant bilinear form on a vertex algebra is
symmetric. 
Our methods, however, are very different from the methods
used in \cite{fhl} and \cite{liform}.

As an application, we introduce a very useful
notion of radical of a vertex algebra. It is equal to  the radical of
a certain canonical invariant bilinear form.  
We prove that a radical-free vertex algebra $A$ is non-negatively
graded, and its component $A_0$ of degree 0 is a commutative
associative algebra, so that all structural maps and operations on $A$
are $A_0$-linear. We also show that in this case 
$A$ is simple if and only if $A_0$ is a field. 

\section{Definitions and notations}\label{sec:vertex}
Here we fix the notations and give some minimal definitions. For more
details on vertex algebras the reader can refer to the books 
\cite{fb,flm,kac2}. 
All spaces and algebras are considered
over a field $\k$ of characteristic 0. We use the following 
notation for divided powers: 
$$
x^{(n)} = 
\begin{cases}
  (x^n)/n!&\text{if}\ n\ge 0\\
0&\text{otherwise}.
\end{cases}
$$ 
\subsection{Definition of vertex algebras}\label{sec:defs}
\begin{Def}\label{dfn:vertex}
A vertex algebra is a linear space $A$ equipped with a family of
bilinear products $a\otimes b \mapsto a(n)b$, indexed by
integer parameter $n$, and with an element $\1\in A$, called the unit,
satisfying the identities (i)--(iv) below.
Let $D:A\to A$ be the map defined by $Da = a(-2)\1$. 
Then \smallskip
\begin{itemize}
\item[(i)] 
$a(n)b = 0$ for $n\gg 0$, \vskip7pt
\item[(ii)]
$\1(n)a = \delta_{n,-1}\,a$ \  \ and \ \ 
$a(n)\1 = D^{(-n-1)} a$,  \vskip8pt
\item[(iii)]
$D(a(n)b) = (Da)(n)b + a(n)(Db)$ \  \ and \ \ 
$(Da)(n)b  = - n\, a (n-1) b$,\vskip8pt
\item[(iv)]
$\displaystyle{a(m)\big(b(n)c\big) - b(n)\big(a(m) c\big) = 
\smash[t]{\sum_{s\ge 0}\binom ms \big(a(s)b\big)(m+n-s)c}}$
\end{itemize}
for all $a,b,c\in A$ and $m,n\in \Z$. 
\end{Def}

This is not the only known definition of vertex algebras. Often the
axioms are formulated in terms of the {\em left adjoint action map}
(a.\thinspace k.\thinspace a. {\em state-field correspondence})
$Y:A\to A[[z,z\inv]]$ 
defined by $Y(a)(z) = \sum_{n\in\Z}a(n)\,z^{-n-1}$, where 
$a(n):A\to A$ is the operator given by   $b\mapsto a(n)b$. 
The most important property of these maps is that they are 
{\em  local}: for any $a,b\in A$ there is $N$ such that 
\begin{equation*}
\big[ Y(a)(w), Y(b)(z)\big]\,(w-z)^N=0.
\end{equation*}

The minimal  $N=N(a,b)$ for which this identity holds is called {\em
the locality index} of $a$ and $b$.
In fact, $N(a,b) = \min\set{n\in\Z}{a(m) b=0\ \forall m\ge n}$.
From (iii) we get  $Y(Da)(z) = \partial_z Y(a)(z)$, which in terms of
coefficients means
\begin{equation}\label{fl:adD}
\ad D{a(m)} = -m\,a(m-1).
\end{equation}

Among other identities that hold in vertex algebras are the associativity
\begin{equation}\label{fl:assoc}
\big(a(n)b\big)(m) c = 
\sum_{s\ge0}(-1)^s \binom ns a(n-s)\big(b(m+s)c\big)
-\sum_{s\le n}(-1)^s \binom n{n-s}
b(m+s)\big(a(n-s)c\big),
\end{equation}
and the quasi-symmetry 
\begin{equation}\label{fl:qs}
a(n)b = - \sum_{i\ge0} (-1)^{n+i} D^{(i)}
\big(b(n+i)a\big).
\end{equation}
For $n\ge0$ the associativity identity \fl{assoc} simplifies to 
\begin{equation*}
\big(a(n)b\big)(m) c = 
\sum_{s=0}^n(-1)^s \binom ns \ad{a(n-s)}{b(m+s)}\,c,
\end{equation*}
which is just another form of the identity (iv) of \dfn{vertex}.

A vertex algebra $A$ is called graded (by the integers) 
if $A = \bigoplus_{d\in\Z}A_d$ is a graded space, 
$A_i(n)A_j\subseteq A_{i+j-n-1}$ and $\1\in A_0$.

It is often required that a vertex algebra $A$ is graded and
$A_2$ contains a special element $\omega$ such that $\omega(0)=D$, 
$\omega(1)|_{A_d} = d$ and the coefficients $\omega(n)$ generate a
representation of the Virasoro Lie algebra:
$$
\ad{\omega(m)}{\omega(n)} = (m-n)\,\omega(m+n-1)+
\delta_{m+n,2}\,\frac 12\binom{m-1}3c 
$$
for some constant $c\in \k$ called the {\em central charge} of $A$. 
In this case $A$ is called {\em conformal vertex algebra} or {\em vertex
operator algebra}, especially when $\dim A_d < \infty$.

\subsection{The action of $sl_2$}\label{sec:sl2}
In order to work with bilinear forms,  we need to deal with vertex
algebras equipped 
with certain additional structure. First of all we will assume that our
vertex algebra $A = \bigoplus_{d\in\Z}A_d$ is graded. We will also need a
locally nilpotent 
operator $D^*:A\to A$  of degree $-1$, such that $D^*\1=0$ and 
\begin{equation}  \label{fl:adDst}
\ad{D^*}{a(m)} =  (2d-m-2)\,a(m+1) + (D^*a)(m)
\end{equation}
for every $a\in A_d$. Let $\delta: A\to A$ be the grading derivation,
defined by $\delta|_{A_d} = d$. It is easy to compute that 
\begin{equation}\label{fl:sl2}
  \ad{D^*}D = 2\delta, \quad
\ad{\delta}D = D,\quad
\ad{\delta}{D^*} =- D^*,
\end{equation}
so that $D^*, D$ and $\delta$ span a copy of $sl_2$. 

An element $a\in A$ such that $D^*a=0$ is called {\em minimal}. 
It is easy to see that if  $A$ is generated by minimal elements, then
any operator $D^*:A\to A$ satisfying \fl{adDst} must be locally
nilpotent.  For any $a\in A$ define $\ord a =\min\set
{k\in \Z_+}{(D^*)^{k+1}a=0}$, so that $\ord a = 0$ for a minimal $a$.
 
Vertex algebras with an action of $sl_2$ as above were called 
{\em quasi-vertex operator algebras} in \cite{fhl} and minimal
elements are sometimes called {\em quasi-primary}. 

If $A$ has a Virasoro element $\omega$, then we always have
$D = \omega(0)$ and $\delta = \omega(1)$, and we can take 
$D^* = \omega(2)$.
 
Throughout this
paper we will always assume that all vertex algebras have an 
$sl_2$-structure, and all ideals, homomorphisms, etc must agree with this
structure. For example, we call a vertex algebra {\em simple} if it
does not have $D^*$-invariant ideals. Note that any $D^*$-invariant
ideal is homogeneous, because it must be stable under the grading
derivation $\delta$.

In the sequel we will need the following easy property of 
$sl_2$-modules of the above type:

\begin{Prop}\label{prop:sl2}\sl
Let $M=\bigoplus_{i\in\Z}M_i$ be a graded module over the Lie algebra 
$sl_2=\k D+ \k\delta + \k D^*$, where $\deg D = 1$, \ $\deg D^* = -1$, \
$D^*$ is locally nilpotent, $\delta|_{M_d} = d$ and the commutation relations
\fl{sl2} hold. Then 
$M_d = (D^*)^{d+1}M_1$ for any $d < 0$.
\end{Prop}
\begin{proof}
Let us first prove that $M_d = D^*M_{d+1}$ for $d<0$.
Take some element $a\in M_d$ of  $\ord a = k$. 
Set 
$$
b = D^*Da - (k+1)(2d-k)\,a = DD^*a- k(2d-k-1)\,a \in M_d.
$$ 
If $a$ is minimal, then $b=0$. 
Otherwise, using that $(D^*)^{k+1}Da = (k+1)(2d-k)\, (D^*)^ka$, we get 
$(D^*)^k\,b=0$, therefore $\ord b < k$. By induction we can assume that
$b\in D^*M_{d+1}$, and hence also $a \in D^*M_{d+1}$, because
$(k+1)(2d-k)\neq 0$ for $d<0$ and $k\ge 0$.

Now the proposition follows from the following claim:
\begin{equation*}
M_0 = D^*M_1 + \ker D^*.
\end{equation*}

Take an element $a\in M_0$ of $\ord a = k$. We will prove by induction
on $k$ that $a \in D^*M_1 + \ker D^*$. If $k=0$ then $a\in \ker
D^*$. Otherwise, take as before $b = D^*Da + k(k+1)a$ so that $\ord b
< k$. By induction, $b\in D^*M_1 + \ker D^*$, and therefore also $a\in
D^*M_1 +\ker D^*$, since $k(k+1)\neq 0$ for $k>0$.
\end{proof}

It follows from \prop{sl2} that $DM_{-1}\subset D^*M_1$. Indeed, if
$a\in M_{-1}$, then $a=D^*b$ for some $b\in M_0$, and we have
$Da = DD^*b = D^*Db$.

For the case when $M=A$ is a vertex algebra  we will obtain that 
$A_d = D^*A_{d+1}$ as a corollary of \lem{Dst} in \sec{form}. 

\begin{Rem}
It is possible (but not easy) to show that if a vertex algebra $A$ is
generated by minimal elements, then  $A_d = D^*A_{d+1}$ for all 
$d\neq 0$.   
In contrast to that, the space
$A_0/D^*A_1$ is non-zero in most interesting cases. We shall see in
\sec{form} that its dual space $\hom(A_0/D^*A_1,\k)$
parameterizes the invariant bilinear forms on $A$.
\end{Rem}

\subsection{The universal enveloping algebra}\label{sec:uea}
For any vertex algebra $A$ we can construct a Lie algebra $L = \cff A$
in the following way \cite{bor,kac2,lihom,freecv}. Consider the linear space
$\k[t,t\inv]\otimes A$, where $t$ is a formal variable. 
Denote $a(n) = a\otimes t^n$ for $n \in \Z$. 
As a linear space, 
$L$ the quotient of $\k[t,t\inv]\otimes A$  by the subspace spanned
by  the  relations $(Da)(n) = -n\,a(n-1)$. The brackets are
given by 
\begin{equation}\label{fl:com}
\ad{a(m)}{b(n)} = \sum_{i\ge0}\binom mi \big(a(i)b\big)(m+n-i),
\end{equation}
which is precisely the identity (iv) of \dfn{vertex}.
The spaces $L_\pm = \spn\set{a(n)}{n
\,\smash{\hbox to0pt{\raisebox{-4pt}{$<$}\hss}\raisebox{4pt}{$\ge$}}\,
0}\subset L$ are Lie subalgebras of $L$ and we have $L = L_-\oplus L_+$.

\begin{Rem}
The construction of $L$ makes use of only the products 
$(n)$ for $n\ge 0$ and the map $D$. This means that it works
for a more generic algebraic structure, known as {\em conformal
algebra} \cite{kac2}.  
\end{Rem}

The  formulas \fl{adD} and \fl{adDst} define derivations $D:L\to L$
and $D^*:L\to L$ so we get an action of $sl_2$ on $L$ by derivations.  
Denote by $\^L = L\ltimes sl_2$ the corresponding semi-direct product.


Now assume that the vertex algebra $A$ is $\Z$-graded. Then so is the
Lie algebra $\^L = \^{\cff A}$ and its universal enveloping algebra
$U = U(\^L)$. The {\em Frenkel-Zhu topology} \cite{fzh} on a homogeneous
component $U_d$ is defined by 
setting the neighborhoods of 0 to be the spaces 
$U^k_d = \sum_{i\le k} U_{d-i} U_i$, so that 
$$
\ldots\subset U_d^{k-1} \subset U_d^k\subset U_d^{k+1}\subset \ldots \subset U_d,
\quad \bigcap_{k\in \Z} U_d^k = 0,
\quad \bigcup_{k\in \Z} U_d^k = U_d.
$$ 
Let 
$\ol U(\^L) = \bigoplus_{d\in\Z} \ol U_d$ be the completion of $U$ in this
topology.  

\begin{Def}\cite{fzh}
The universal enveloping algebra $U(A)$ of a graded 
vertex algebra $A$ is the
quotient of $\ol U(\^{\cff A})$ modulo the ideal generated by the
relations
$$
\big(a(n)b\big)(m) = 
\sum_{s\ge0}(-1)^s \binom ns a(n-s) b(m+s)
-\sum_{s\le n}(-1)^s \binom n{n-s}
b(m+s)a(n-s)
$$  
for $a,b\in A$ and $m,n\in \Z$.
\end{Def}
Note that the relations above are simply the associativity identity
\fl{assoc}. It is enough to impose these relations only for $n<0$,
because for $n\ge 0$ they are  equivalent to
\fl{com}.

A graded space $M$ is a module over a vertex algebra $A$ if it is a
$U(A)$-module that agrees with the Frenkel-Zhu topology on $U(A)$,
i.e. a convergent sequence of elements from $U(A)$ makes sense as an
operator on $M$. 
\section{The involution on the universal enveloping algebra}
\label{sec:involution}
Let $A$ be a graded vertex algebra with $sl_2$-structure, as in \sec{sl2}.

\begin{Prop}\mylabel{prop:adjoint}\sl
Set
\begin{equation}\label{fl:adjoint}
 a(m)^* = (-1)^{\deg a} \sum_{i\ge 0} 
\big(D^{*(i)}a\big)(2\deg a-m-2-i)
\end{equation}
for a homogeneous $a\in A$ and $m\in \Z$, and also $D^{**}=D$ and
$\delta^*=\delta$. Then there is a continuous anti-involution 
on the enveloping algebra $U(A)$ defined by  $x \mapsto x^*$. 
\end{Prop}
Note that $U(A)_d^*\subseteq U(A)_{-d}$. 

\begin{Rem}
Let $M$ be a module over $A$. Assume that \prop{adjoint} is true. Then
the graded dual space $M' = \bigoplus_{d\in\Z}\hom(M_d,\k)$ also becomes
an $A$-module, called {\em the contragredient module} of $M$, where
the action of $A$ is given by $\big(a(n)f\big)(v) = f\big(a(n)^*v\big)$,
\ $(Df)(v) = f(D^*v)$ and $(D^*f)(v) = f(Dv)$ for $a\in A$, \ $f\in
M'$, \ $v\in M$ and $n\in \Z$. Contragredient
modules were introduced in \cite{fhl} without the use of
\prop{adjoint}. In fact, it is easy to show that \prop{adjoint} is
equivalent to the existence of the contragredient
module. For the sake of completeness though, we present an independent
proof of \prop{adjoint}, that avoids
the use of formal calculus.
\end{Rem}

\begin{proof}[Proof of \prop{adjoint}]
It is straightforward to check that  $a(m)^{**} = a(m)$. This is
particularly easy to see when $D^*a=0$, and that already settles the case
when $A$ is generated by minimal elements. The general case is not
much more difficult. 

Let us show that \fl{adjoint} defines an anti-endomorphism of
$U(A)$. We do this in three steps.

\subsubsection*{Step 1} We show first that \fl{adjoint} defines an
anti-endomorphism of the coefficient Lie algebra 
$\^L = \cff(A)\ltimes sl_2$. 
We have to check that  
\begin{equation}\label{fl:tocheck}
\begin{gathered}
\ad{b(n)^*}{a(m)^*} = \sum_{s\ge 0} \binom ms 
\big(a(s)b\big)(m+n-s)^*,\\
\ad{a(m)^*}{D^*} = -m\,a(m-1)^*
\end{gathered}
\end{equation}
for  $a,b$ running over a set of generators of $A$ and $m,n\in \Z$.
Let us check this  when $a$ and $b$ are minimal, thus  proving
\fl{tocheck} for the case when $A$ is generated by the minimal
elements.  The general case is not much different.

Let $\deg a = d,\, \deg b = e$ and  $D^*a = D^*b = 0$. Then we get
\begin{align*}
\ad{b(n)^*}{a(m)^*} &= -\ad{(-1)^d a(2d-m-2)}{(-1)^e b(2e-n-2)} \\
&=(-1)^{d+e+1}\sum_{s\ge 0} \binom{2d-m-2}{s}  
\big(a(s)b\big)(2d+2e-m-n-s-4)  
\end{align*}
On the other hand, using that 
$$
D^{*(i)}\big(a(s)b\big) = \binom{2d-s-2}i\
a(s+i)b,
$$
we have
\begin{align*}
\sum_{s\ge 0} &\binom ms \big(a(s)b\big)(m+n-s)^* \\
&=  \sum _{i,s\ge 0} (-1)^{d+e+s+1} \binom ms\binom{2d-s-2}i
\big(a(s+i)b\big)(2d+2e-m-n-i-s-4)\\
&= (-1)^{d+e+1}\sum_{j\ge0}
\sum_{s=0}^j (-1)^s \binom ms\binom{2d-2-s}{j-s}
\big(a(j)b\big)(2d+2e-m-n-j-4),
\end{align*}
and the first equality of \fl{tocheck} follows from a binomial
identity. Similarly, by \fl{adDst} and \fl{adjoint},
$$
\ad{D^*}{a(m)^*} = (-1)^d\, \ad{D^*}{a(2d-m-2)}
= (-1)^d\, m \, a(2d-m-1) = m\,a(m-1)^*,
$$
and this proves the second equality of \fl{tocheck}.

\subsubsection*{Step 2}
Assume that we have two vertex algebras $A$ and $B$, both with 
$sl_2$-structures as above. 
Assume also that for $B$ the proposition is known, that is, formula
\fl{adjoint} defines an anti-involution of $U(B)$, and suppose 
we have a surjective homomorphism $\f:B\to A$ of vertex algebras,
which preserves the $sl_2$-structure. Then the proposition also holds
for $A$. 

Indeed, we can extend $\f$ to a homomorphism $U(\f):U(B)\to U(A)$ of
enveloping algebras. Since $\f$ is surjective, so is $U(\f)$, and
therefore $U(A) = U(B)/\ker U(\f)$. It is easy to see that $\ker
U(\f)$ is the ideal of $U(B)$ generated by the coefficients $a(m)$ for
all $a\in \ker \f \subset B$ and $m\in \Z$. But if $a\in \ker \f$,
then also $D^*a\in \ker \f$, therefore $a(n)^* \in \ker U(\f)$, and it
follows that $\big(\ker U(\f)\big)^* \subseteq \ker U(\f)$.

\subsubsection*{Step 3}
Consider the Verma $\^L$-module 
$B = U(\^L)\otimes_{U(L_+)\ltimes sl_2}\k \1$ generated by 
a vector $\1$ such that $L_+\1=0$ and also $D\1=D^*\1=0$.  Then $B$
has a structure of vertex algebra such that $U(B)$ is the
completion in the Frenkel-Zhu topology of the universal enveloping
algebra $U=U(\^L)$ of $\^L$. The action of $D^*$ on $B$ is given by
\fl{adDst}, therefore Step 1 implies that \fl{adjoint} defines an
anti-involution of $U$. 

Now we show that this anti-involution on $U$ is continuous in the Frenkel-Zhu
topology, and therefore it can be extended to $U(B) = \ol U$. 
Indeed, let $u_1,u_2, \ldots \in U_d$ be a convergent
sequence. This means that for any $k\in \Z$ all elements $u_n$ with
sufficiently large $n$ belong to 
$U_{d-i}U_i$ for some $i\le k+d$. But then
$u_n^*\in \big(U_{d-i}U_i\big)^*\subseteq U_{-i} U_{i-d}\subseteq U_{-d}^k$,
therefore the sequence $\{u_n^*\}$ also converges in the Frenkel-Zhu
topology. 

 Since $A$ is also an $\^L$-module
generated by an element $\1$ such that $L_+\1=D\1=D^*\1=0$, there
is a unique $\^L$-module homomorphism $\f:B\to A$ such that $\f(\1)=\1$. It
is well known \cite{dlm2001,primc,freecv} that $\f$ must be a vertex
algebra homomorphism preserving the $sl_2$-module
structures. Therefore, Step 2 implies that \fl{adjoint} gives an
involution of $U(A)$.
\end{proof}

\section{Invariant bilinear forms}\label{sec:form}
Let as before $A$ be a vertex algebra with $sl_2$-structure. 
Let $V$ be a vector space over  $\k$.
\begin{Def}\label{dfn:ibf}
A $V$-valued bilinear form $\form\cdot\cdot$ on $A$ is called {\it invariant} if
$$\form {a(m)b}c = \form b{a(m)^*c}\quad\text{and}\quad 
\form {Da}b = \form a{D^*b}
$$  
for all $a,b,c\in A$ and $m\in\Z$.
\end{Def}

The radical 
$\rad\form\cdot\cdot = \set{a\in A}{\form ab=0\ \forall b\in A}$
of an invariant form is an ideal of $A$. Also, we have 
$$
\form{\delta a}b = \frac 12 \bform{\ad{D^*}{D}a}b 
=  \frac 12 \bform{a}{\ad{D^*}{D}b} = \form a{\delta b},
$$
and this implies that $\form {A_i}{A_j}=0$ for $i\neq j$.

When $A$ is a highest weight module over affine or Virasoro Lie
algebra $\goth g$ (see e.g. \cite{fzh,kac2}), then one can choose
$D^*$ so that the involution
$^*$ of \sec{involution} is the extension of the Chevalley
involution of $\goth g$, and the 
contravariant form on $A$ (see e.g. \cite{kac1}) 
will be invariant in the sense of \dfn{ibf}.
See also the example in \sec{heisen}.

\begin{Rem}
  Assume that $A$ has an  invariant bilinear form
$\formdd$. Viewing $A$ as a module over itself, consider the
contragredient module $A'$, see the Remark following \prop{adjoint} in
\sec{involution}. Then there is an $A$-algebra homomorphism 
$\nu:A\to A'$ given by  $a\mapsto  \form a\cdot$. In particular, if
$\form \cdot\cdot$ is  non-degenerate,
and all the homogeneous components of
$A$ are finite dimensional, then $\nu:A\to A'$ is an isomorphism. This
is a key observation in \cite{liform}.
\end{Rem}

Given an invariant form $\form \cdot\cdot$ on $A$, one can consider a
linear map $f: A_0\to V$ defined by $f(a)= \form\1 a$. Clearly,
$f(D^*A_1)=0$, so that $f\in \hom (A_0/D^*A_1, V)$, and $f$
defines the form uniquely. 
We show that this in fact gives an isomorphism of the space of all
invariant bilinear forms with $\hom (A_0/D^*A_1, V)$.

\begin{Thm}\mylabel{thm:sym}\sl
Every invariant bilinear form on $A$ is symmetric and there is a
one-to-one correspondence between linear maps
$f: A_0/D^*A_1 \to V$ and $V$-valued invariant bilinear forms
$\form\cdot\cdot$ on $A$, given by $f(a)= \form\1 a$ for $a\in A_0$.
\end{Thm}

For the proof we need the following lemma.

\begin{Lem}\mylabel{lem:Dst}\sl
For any $a,b \in A$, if either $m\ge 0$ or $m<-\ord b-1$, then 
$a(m)^* b \in D^*A$. 
\end{Lem}

\begin{proof}
Set $\deg a = d$ and $\ord b = k$. Denote $\big((D^*)^ja\big)(l) =
g_{2(d-j)-l-1}^{d-j-l-1}\in U(A)$, so that 
$$
g_m^n =
\big((D^*)^{d+n-m}a\big)(m-2n-1).
$$
For an operator $u\in U(A)$ write $u\sim 0$ if $ub \in  D^*A$.
Then by \fl{adDst} we have
$$
m\,g_m^n + g_{m-1}^n + g_{m+1}^{n+1}D^*\sim 0.
$$
Consider the generating function $g(x,y) = \sum_{m,n\in \Z} g_m^n\,x^m
y^n$. Let $R:U(A)\to U(A)$ be the operator given by $Ru =
uD^*$. Note that $R^i \sim 0$ for $i > k$, since $(D^*)^{k+1}b = 0$.
Then the above relation reads as
$$
x\,\partial_x g + x g+\frac R{xy} \,g \sim 0,
$$
therefore
$$
g(x,y) \sim  \exp\(-x+\frac R{xy}\) \, g_0(y)  
$$
for some series $g_0(y)\in U(A)[[y^{\pm1}]]$.
On the other hand, 
$$
a(m)^* = (-1)^d \sum_{i\ge 0}  \big(D^{*(i)}a\big)(2d-m-2-i)
= (-1)^d \sum_{i\ge 0} \frac 1{i!}\, g_{m+1-i}^{m+1-d}.
$$
Set $h_m^n = \sum_{i\ge 0} \frac 1{i!}\, g_{m-i}^n$, so that 
$a(m)^* = (-1)^d\, h_{m+1}^{m+1-d}$. The generating function for
$h_m^n$'s is 
$$
h(x,y) = \sum_{m,n\in\Z} h_m^n\,x^my^n = e^x g(x,y) \sim \exp
\(\frac R{xy}\) \, g_0(y).
$$ 
Since $(R)^i \sim 0$ for $i > k$, the only powers of $x$ with
coefficients $\not\sim0$ in $h(x,y)$ are $x^0$, $x^{-1}$, $\ldots$,
$x^{-k}$, therefore $h_m^n \sim 0$ if either $m>0$ or $m<-k$, and the
lemma follows.
\end{proof}

\begin{Rem}
As a corollary, we get a proof of the fact that $A_d = D^*A_{d+1}$
different from that of\prop{sl2}.  Indeed, take an element $a\in A_d,
\ d < 0$. Write $a = a(-1)\1$. Then $a(-1)$ is the dual of $\pm
\sum_{i\ge 0}(D^{*(i)}a)(2d-i-1)$. But $\ord \1 = 0$ and $2d-i-1 \neq
-1$ for all $i\ge 0$, so \lem{Dst} implies that $a(-1)\1\in
D^*A_{d+1}$.
\end{Rem}

\begin{proof}[Proof of \thm{sym}]
We need to show that for every linear map $f:A_0/D^*A_1\to V$
there is a symmetric invariant bilinear form $\form\cdot\cdot$ on $A$ such
that $f(a) = \form \1a$ for all $a\in A_0$.

For homogeneous $a,b\in A$, define $\form ab = f\big(a(-1)^*b\big)$
if $\deg a = \deg b$ and $\form ab = 0$ otherwise.
Let us show first that $\form {a(n)b}c = \form b{a(n)^*c}$
for all $a,b,c\in A$ and $n\in \Z$. 
Indeed, it is enough to check that 
$\big(a(n)b\big)(-1)^*c \equiv b(-1)^*a(n)^*c \mod D^*A$.
But using the associativity \fl{assoc}, we have 
$\big(a(n)b\big)(-1) - a(n)b(-1)\in U(A)L_+$, 
and by \lem{Dst}, $u^*A\subseteq
D^*A$ for any $u\in U(A)L_+$. Similarly, $\form{Da}b = f(a(-1)^*D^*b)
= \form a{D^*b}$. This proves that the form is invariant.

Now we show that the form $\form\cdot\cdot$ is symmetric. We need to
have $a(-1)^*b \equiv b(-1)^*a \mod D^*A$ for every $a,b\in A$ of the
same degree. Set
$u = a(-1)^*b(-1) \in U(A)$ so that $\deg u = 0$. We have to show that
\begin{equation}\label{fl:uust}
(u-u^*)\1 \in D^*A_1.
\end{equation}
By \fl{adjoint}, $u$ is a linear combination of operators of
the form $c(m)b(-1)$ for some $c\in A$ and $m\in \Z$. But 
$c(m)b(-1) \equiv \big(c(m)b\big)(-1) \mod
U(A)L_+$. 
Hence it is enough to prove \fl{uust} for the cases when $u\in U(A)L_+$
and when $u=a(-1)$ for some $a\in A_0$. 

If $u\in U(A)L_+$ then $u\1=0$ and $u^*\1\in D^*A$ by \lem{Dst}.
If $u=a(-1)$ for  $a \in A_0$, then 
$$
a(-1)-a(-1)^* = 
\big(a(-1)^*-a(-1)\big)^* 
= \sum_{i\ge 1}\big(D^{*(i)}a\big)(-i-1)^*
$$
and again by \lem{Dst} we conclude that
$\big(a(-1)-a(-1)^*\big)\1\in D^*A_1$, since $\ord \1 = 0$.
\end{proof}

\begin{Rem}
It follows from \thm{sym} that if $u\in U(A), \ \deg u = d,$ is such
that $u\1 = 0$, 
then $u^*A_d\in D^*A_1$. Indeed, $0=\form {u\1}{A_d} = \form\1{u^*A_d}$ for
every invariant bilinear form on $A$. In fact, using techniques from
\cite{freecv,cfva} one can show that $\set{u\in U(A)}{u\1=0} = U(A)L_+$.
\end{Rem}

\section{Radical of a vertex algebra}\label{sec:radical}
\subsection{The ideal generated by $D^*A_1$}\label{sec:I}
Let $A$ be a vertex algebra with an $sl_2$-structure as before. 
Let $I\subset A$ be the ideal generated by $D^*A_1$. Denote as usual 
$I_d = I\cap A_d$. 


\begin{Prop}\label{prop:I}\sl
\begin{enumerate}
\item\label{I:neg}
$I_d = A_d$ for $d<0$.
\item \label{I:0}
$I_0 = A(-1)D^*A_1$.
\item\label{I:asscom}
$A_0/I_0$ is an associative commutative algebra with respect to the
operation $a\otimes b\mapsto ab=a(-1)b$, and $\1$ is its unit. 
\end{enumerate}
\end{Prop}

Denote $E =  A(-1)D^*A_1$. 
For the proof of \prop{I} we need the following lemma.

\begin{Lem}\label{lem:I}\sl
Let $a\in A_i$ and $b\in A_j$ for $i,j\le 0$, and if $j=0$, then $b\in
D^*A_1$. Then $a(i+j-1)b \in E$.  
\end{Lem}
\begin{proof}
By \prop{sl2} we have $a = (D^*)^i a_0$ and $b = (D^*)^j b_0$ for some 
$a_0\in A_0$ and $b_0\in D^*A_1$. Applying \fl{adDst} $-i$ times we get that 
$$
a(i+j-1)b \equiv (-1)^i a_0(i+j-1)(D^*)^{-i-j}b_0 \mod E.
$$

If $d=i+j=0$, then $a(i+j-1)b=a(-1)b\in E$ by the definition. 
We will show by induction on $d$ that
$a(d-1)(D^*)^{-d}b \in E$ for any $a\in A_0$, \ $b\in D^*A_1$ and $d<0$.

So assume that $a(k-1)(D^*)^{-k} b\in E$ for all $0\ge k>d$.
Then, using \fl{adD}, we get 
$$
E\supset D^*A_1\ni Da(d)(D^*)^{-d}b = a(d)D(D^*)^{-d}b - d\, a(d-1)(D^*)^{-d} b. 
$$
By \prop{sl2}, $D(D^*)^{-d}b= (D^*)^{1-d}b_1$ for some $b_1\in D^*A_1$,
so by induction, $a(d)D(D^*)^{-d} b \in E$, and hence also
$a(d-1)(D^*)^{-d} b\in E$.
\end{proof}

\begin{Rem}
In fact, \lem{I}, combined with \lem{Dst}, shows that $E = U(A)_0D^*A_1$.
\end{Rem}

\begin{proof}[Proof of \prop{I}]
(\ref{I:neg})\ 
follows immediately from \prop{sl2}.

\smallskip\noindent
(\ref{I:0})\quad 
Let $f:A_0\to A_0/E$ be the canonical projection, and let 
$R = \rad\formdd$ be the radical of the corresponding form on $A$. 
We claim that $R_0=R\cap A_0 = E$. Indeed, if $a\in A_0\ssm E$, then 
$f(a)=\form a\1\neq 0$ and therefore $a\not\in R$. In the other
direction, take $a\in D^*A_1$. Then for an arbitrary $b\in A_0$ we
have 
$$
\form ab = f\big(b(-1)^*a\big) = 
\sum_{i\ge 0} f\big((D^{*(i)}b)(-1-i)a\big) =0,
$$
since $(D^{*(i)}b)(-1-i)a \in E$ by \lem{I}. Therefore, $D^*A_1
\subset R$, and since $R$ is an ideal, we must have $E\subset R$. 
It follows that $I\subset R_0 = E$, on the other hand obviously $E\in
I$, and that proves (\ref{I:0}).

\smallskip\noindent
(\ref{I:asscom})\quad 
Using \fl{assoc} and \lem{I} we get
$$
\big(a(-1)b\big)(-1)c = \sum_{s\ge 0} a(-1-s)b(-1+s)c + 
\sum_{s\ge 0} b(-2-s)a(s)c \equiv a(-1)b(-1)c \mod I_0
$$ 
for any $a,b,c\in A_0$. By \fl{qs} we have 
$a(-1)b \equiv b(-1)a\mod DA_{-1}$, and we know that $DA_{-1}\subset D^*A_{i+1}
\subset I_0=E$. So the statement (\ref{I:asscom}) follows.
\end{proof}

\subsection{The radical}\label{sec:rad}

\begin{Def}\label{dfn:radical}
The radical $\rad A$ of a vertex algebra $A$ is the radical
$\rad\formdd$ of the invariant bilinear form $\formdd$ corresponding
to the projection $f:A_0\to A_0/A_0(-1)D^*A_1$.
\end{Def}

\begin{Rem}
This definition has nothing to do with the radical defined in \cite{dlmm}.  
\end{Rem}

Recall that by \propp{I}{0}, the space $I_0=A_0(-1)D^*A_1$ is the 
degree 0 component of the ideal $I\subset A$ generated by $D^*A_1$. From
the proof of \propp{I}{0} it follows that $(\rad A)_0 = I_0$. 
Moreover, the following  is true:
\begin{Prop}\label{prop:rad}\sl
\begin{enumerate}
\item\label{rad:max} 
The radical $R = \rad A$ is the maximal among all the
ideals $J\subset A$ such that $J_0 = A_0(-1)D^*A_1$. 
\item\label{rad:rad} 
$\rad(A/R)=0$.
\end{enumerate}
\end{Prop}

\begin{proof}
Let $J\subset A$ be an ideal such that $J_0 =  A_0(-1)D^*A_1$. 
Then for any $a\in J_d$ and $b\in A_d$ we get 
$\form ab = f\big(b(-1)^*a\big) = 0$, therefore $a\in R$. This
proves (\ref{rad:max}). The statement (\ref{rad:rad}) follows from the
fact that the form $\formdd$ becomes non-degenerate on $A/R$.
\end{proof}

\subsection{Radical-free algebras}\label{sec:radfree}
Assume that $R =\rad A \subsetneqq A$. We would like to investigate the
quotient vertex algebra 
$A/R$. By \propp{rad}{rad} we have 
$\rad A/R = 0$,  so we can assume that $\rad A =0$. Then $D^*A_1 = I =
0$, the projection $f:A_0 \to A_0/I_0$  is the
identity map, and the corresponding  canonical  invariant
$A_0$-valued bilinear form $\formdd$ is non-degenerate.

The following theorem summarizes the properties of
radical-free algebras.

\begin{Thm}\label{thm:rad0}\sl
Let $A$ be a vertex algebra such that $\rad A = 0$.
\begin{enumerate}
\item\label{rad0:neg} 
$A_d=0$ for $d<0$.
\item\label{rad0:0}
$A_0$ is  an associative
commutative unital algebra with respect to the product $a\otimes b \mapsto ab
= a(-1)b$, with unit $\1$.
\item\label{rad0:mod}
$A$ is a vertex algebra over the associative
commutative algebra $A_0$. 
\item\label{rad0:ideals}
Let $I\subset A$ be an ideal. Then $I\cap A_0 = 0$ if and only if
$I=0$. 
\item\label{rad0:simple}
For any ideal $J_0\subset A_0$ the space $J_0A$ is an ideal of $A$
such that $J_0A\cap A_0 =J_0$. It follows that   
$A$ is simple if and only if $A_0$ is a field. 
\end{enumerate}
\end{Thm}

Note that all the definitions in \sec{vertex} make sense when $\k$
is a commutative ring containing $\Q$. So the statement
(\ref{rad0:mod}) means that elements from $A_0\subset A$ behave like
constants, i.e. $A_0$ acts on each component $A_d$ by $ab=a(-1)b$ for
$a\in A_0$, \ $b\in A_d$,  all
structural maps $a(m), D, D^*:A\to A$ are $A_0$-linear,  and therefore 
$DA_0 = 0$ and $A_0(n)A_d=0$ for $n\neq -1$. Also, the canonical
non-degenerate $A_0$-valued form
$\formdd$ on $A$ is $A_0$-bilinear.

Recall that a vertex algebra with $sl_2$-structure is called simple if
it does not have any $D^*$-invariant ideals. 

\begin{proof}
The statements (\ref{rad0:neg}) and (\ref{rad0:0}) follow from
\propp{I}{neg} and (\ref{I:asscom}) respectively, while 
(\ref{rad0:ideals}) follows from \propp{rad}{max}.

Since $D^*A_1=0$, we have $\form{Da}b =0$ for any  and
$a\in A_0$,\ $b\in A_1$, therefore $DA_0\subset \rad A = 0$. But then
$a(n)=0$ for $n\neq -1$ as an element in the Lie algebra $\cff A$, see
\sec{uea}. Now the statement (\ref{rad0:mod}) follows easily from the
formulas \fl{adD}, \fl{assoc} and \fl{adDst} in \sec{defs}.

Finally, $I_0A$ is an ideal of $A$ by (\ref{rad0:mod}), and therefore,
if $A$ is simple, then $A_0$ must be a field. This proves one
direction of (\ref{rad0:simple}). The other direction follows from
(\ref{rad0:ideals}).
\end{proof}


\begin{Rem}
Let $A$ be a vertex algebra with $\rad A=0$ as above. Since $DA_0=0$,
the formulas \fl{assoc} and \fl{qs} imply that $A_1$ is a
Lie algebra with the bracket  $\ad ab = a(0)b$. The
bilinear form $\formdd$ on $A$  is an invariant
symmetric $A_0$-bilinear form on that Lie algebra, and we have 
$\form ab = a(1)b$ for any $a,b\in A_1$. The   Lie algebra $A_1$ acts
on  $A$ by derivations of degree 0, so that $ab = a(0)b$ for $a\in A_1$
and $b\in A$. If  $A$ is simple, then 
$A_1(0)A=0$, since  $A_1(0)A$ is always an ideal of $A$ and $\1\not\in
A_1(0)A$.  
\end{Rem}




If $\k=\C$ and $A$ is a simple radical-free vertex algebra, generated
by countably many elements, then we must have $A_0=\C$. 
We see that radical-free vertex algebras are very close to be of a CFT
type, after we set $\k=A_0$.  The only
thing that might be missing  is the condition $\dim A_d<\infty$.

\medskip\noindent
\bf Question.\rm\quad
Are there any radical-free vertex algebras $A$ with
$sl_2$-structure, such that the homogeneous components $A_d$ have
infinite rank over $A_0$?

\section{Examples}
\subsection{Heisenberg algebra}\label{sec:heisen}
Let $A$ be the Heisenberg vertex algebra. By definition, it is
generated by a single element $a\in A_1$ such that the locality 
$N(a,a)=2$ and $a(0)a = 0$, \ $a(1)a = \1$. One can
show (see e.g. \cite{cfva}) that these conditions define $A$
uniquely. We have $A_0 = \k\1$ and $A_1 = \k a$.  It is well known,
that $A$ is a simple vertex algebra in the sense that it does not have any ideals
at all, regardless of any $sl_2$-structure. 

For every $k\in \k$ the element 
$\omega_k = \frac 12\, a(-1)a + k\,Da\in A_2$ 
is a Virasoro element of $A$. If we set $D^* = \omega_0(2)$, 
then $D^*a=0$, hence $A_0/D^*A_1 = \k\1$, and 
therefore $A$ has only one invariant
bilinear form $\formdd$, up to a scaling factor, and in this case 
$\rad A = \rad\formdd=0$. 
If instead we set 
$D^* = \omega_k(2)$ for $k\neq 0$, then  $D^*a=-2k\1$. In
this case $A_0/D^*A_1 = 0$, \ $\rad A = A$, and $A$ does not have any invariant
bilinear forms.

\subsection{Free vertex algebras}\label{sec:free}
Let $\cal G$ be a set and $N:\cal G\times \cal G\to \Z$ be a 
symmetric function such that $N(a,a)\in 2\Z$. 
In \cite{cfva} we have constructed a vertex
algebra $F = F_N(\cal G)$ which we called {\em 
the free vertex algebra generated by $\cal G$ with respect to locality
bound $N$}.  It is generated by $\cal G$ so that the locality 
of any $a,b\in \cal G$ is exactly $N(a,b)$, and for any other
vertex algebra $A$ generated by $\cal G$ with the same bound on
locality, there is a homomorphism $F\to A$ preserving $\cal G$. The
algebra $F$ is graded by the semilattice $\Z_+[\cal G]$, we refer to
this grading as the grading by weights. Also, $F$ is graded by degrees
if we set $\deg a = -\frac 12 N(a,a)$ for $a\in \cal G$. If we denote 
by $F_{\lambda, d}$ the space of all elements of $F$ of weight
$\lambda \in \Z_+[\cal G]$ and degree $d\in \Z$, then 
$$
F_{\lambda, k}(n)F_{\mu,l} \subseteq
F_{\lambda+\mu,k+l-n-1},\quad DF_{\lambda, k} \subseteq
F_{\lambda,k-1}\quad\text{and}\quad \1\in F_{0,0}.
$$

It is not difficult to see that we can define an operator $D^*:F\to
F$, satisfying \fl{adDst}, by setting $D^*a=0$ for all $a\in \cal G$. 
Then we have $D^*F_{\lambda, d} \subseteq
F_{\lambda,d+1}$. 

It is shown in \cite{cfva} that for any weight
$\lambda = a_1+\ldots +a_l\in \Z_+[\cal G]$
there is the minimal degree 
$$
d_{\min}(\lambda) =-\frac12 \sum_{i,j=1}^lN(a_i,a_j)\in \Z
$$ 
such that 
$\dim_\k F_{\lambda, d_{\min}(\lambda)} = 1$ and    
$\dim_\k F_{\lambda, d} = 0$ for $d <  d_{\min}(\lambda)$.
For $d > d_{\min}(\lambda)$ the dimension of $F_{\lambda, d}$ 
is equal to the number of partitions of $d-d_{\min}(\lambda)$ into a
sum of $l$ non-negative integers colored by $a_1, \ldots, a_l$.

It is also proved in \cite{cfva} that $F$ is embedded into the vertex
algebra $V_\Lambda$ corresponding to the lattice 
$\Lambda = \Z[\cal G]$. The scalar product on
$\Lambda$ is defined by $\form ab = -N(a,b)$ for $a,b\in\cal G$.

Since $D^*$ is homogeneous with respect to both weights and degrees,
and the only element of weight 0 in $F$ is $\1$, up to a scalar, it
follows that $\1 \not\in D^*F_1$. So one can always define a
functional $f:F_0/D^*F_1\to \k$ so that $f(\1)=1$ and 
$f(F_{\lambda,0}) = 0$ for every weight $\lambda \neq 0$. The corresponding
invariant bilinear form on $F$ is a restriction of  the canonical
bilinear form on $V_\Lambda$.
If  $N(\cal G,\cal G)< 0$, then this will be the only form on $F$, up
to a scalar. If, on the contrary, the
localities $N(a_i,a_j)$ are large, then $\dim F_{\lambda,0}$ grows as 
$\lambda$ gets longer, and  it follows from the results of \cite{cfva}
that $\dim F_{\lambda,0}/(\rad F)_{\lambda,0}>0$ when
$\lambda\in \Z_+[\cal G]$ is big enough.
Hence, in this case there are infinitely many different invariant
bilinear forms on 
$F$, and infinitely many simple radical-free vertex
algebras, which are homomorphic images of $F$. 

\begin{Rem}
It appears that the algebra $\bar F = F/\rad F$ is a very interesting
object from combinatorial point of view. It is possible to show that 
the commutative associative algebra $\bar F_0$ is isomorphic to a
polynomial algebra in infinitely many variables. Also, the Question at
the end of \sec{radical} is equivalent to asking whether the components 
$\bar F_d$ have finite rank as modules over $\bar F_0$. 
\end{Rem}

\subsection*{Acknowledgments}
I thank Haisheng Li for inspiring discussions and for pointing out an
inaccuracy in an earlier version of this paper.

\bibliography{../vertex,../my,../general,../conformal}

\end{document}